\documentclass[12pt]{article}
\usepackage{graphicx,color}

\usepackage{anysize,url}
\usepackage{amsfonts}
\usepackage{amsmath}
\title{An Upper Bound for the\\Number of Planar Lattice Triangulations}
\author{Emile E. Anclin\\
Institute of Mathematics, MA 6-2\\
TU Berlin,
D-10623 Berlin, Germany\\
\url{anclin@math.tu-berlin.de}}

\newtheorem{thm}{Theorem}

\newtheorem{lemma}[thm]{Lemma}

\newtheorem{defn}[thm]{Definition}

\def\endofpr{$\square$\\}

\def\be{\begin{eqnarray}}
\def\ee{\end{eqnarray}}

\def\bee{\begin{eqnarray*}}
\def\eee{\end{eqnarray*}}

\def\raw{\rightarrow}

\def\Z{\mathbb{Z}}

\def\int{\mathop{\rm int}}
\def\area{\mathop{\rm area}}
\def\conv{\mathop{\rm conv}}

\begin{document}

\maketitle

\begin{abstract}%
We prove an exponential upper bound for the number $f(m,n)$ of all
maximal triangulations of the $m\times n$ grid:
\[
f(m,n)\ \ <\ \  2^{3mn}.
\]
In particular, this improves a result of S. Yu. Orevkov \cite{Ore}.
\end{abstract}

We consider  lattice polygons $P$ (with vertices in $\Z^2$), for
example the convex hull of 
the grid $P_{m,n}:=\{ 0,1,\dots ,m\}\times\{0,1,\dots ,n\}$.
We want to estimate the number of \textbf{maximal} lattice triangulations of
$P$, i.e., triangulations  using all integer points $P\cap\Z^2$
in~$P$. These are exactly the 
\textbf{unimodular} triangulations, in which all the triangles have
integer vertices and area $\frac{1}{2}$.
From now on we will talk only about unimodular  triangulations.
Denote by  $f(P)$ the number of (unimodular) triangulations of $P$ and
by $f(n,m)$ the number of  triangulations of $P_{m,n}$. 
S. Yu. Orevkov's  upper bound \cite{Ore} 
is $f(m,n)\leq 4^{3mn}$.

\begin{thm}\label{thm:hay}
The number $f(P)$ of maximal triangulations of a lattice polygon $P$ is
bounded by\vskip-10mm 
\[ f(P)\ \ \le\ \ 2^{|E'|},\]
where $|E'|$ is the number of inner (non-boundary) edges in \emph{any}
unimodular triangulation of~$P$.

In particular, the number of unimodular triangulations of the grid $P_{m,n}$ is bounded by
\[ f(m,n)\ \ \le\ \ 2^{3mn-m-n}\ <\ 2^{3mn}.\]
\end{thm}

\section*{The Haystack Approach} 

Let $P$ be a closed, not necessarily convex lattice polygon and $\int(P)$ its interior.
 Define $M:=(\frac{1}{2} \Z^2{\setminus}\Z^2)\,\cap\, \int(P)$, the possible 
midpoints of the inner edges of a lattice triangulation of $P$.

\begin{lemma}
In any unimodular triangulation of  $P$  there is a canonical
bijection between the inner edges $E'$ and their midpoints in $M$.
\end{lemma}

\textsc{Proof}: The injection from $E'$ to $M$ is clear. \\
On the other hand all unimodular  triangles are $SL(2,\Z)$-equivalent to 
$\Z^2$-translates of
$\conv\{{\bf 0},e_1,e_2\}$, so they don't contain interior points
from~$M$.  
\endofpr

\textsc{Notation}:
For  a subcomplex $S$ of a  triangulation of $P$ and $r\in M $, 
  if there is an edge through $r$ in $S$ we denote it by $e_S(r)$.  We use a lexicographic order on $(\frac{1}{2}\Z)^2$:    
\[
 (x_1,y_1) \prec (x_2,y_2) \ \ \Longleftrightarrow\ \ 
 [y_1<y_2] \mbox{ or } [ y_1= y_2\mbox{ and }x_1<x_2].  
\]

\begin{defn}
A \textbf{haystack} $H$ (with respect to some $r\in M$) is a
subcomplex  of a triangulation of $P$, with  an edge through 
$r'\in M$  if and only if $r'\prec r$. 
\end{defn}

\textsc{Proof of Theorem~\ref{thm:hay}:}\\
The idea is to run through $M$  lexicographically, and at each step to  add an edge through $r\in M$.
 We will see that in each step there are at most two possibilities to put  the new  edge through $r$.
\begin{figure}[ht]
\begin{center}
\input{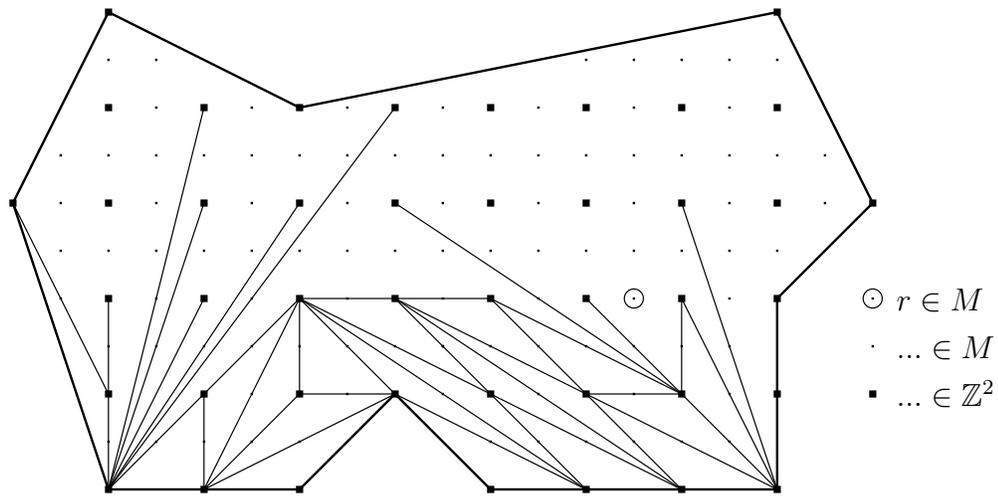}
\caption{A haystack with respect to $r$.}
\end{center}
\end{figure}

We proceed by induction on the totally ordered set $(M,\prec)$, thus
proving that the number of haystacks  with respect to some $r\in M$ is
$\leq 2^{e_r}$ where $e_r$ is the number of predecessors  of $r$ in~$M$.  
Thus after the final step (that is, after processing the largest
$r$ in $(M,\prec)$) we have obtained that there are at most
$2^{|M|}=2^{|E'|}$ unimodular triangulations of~$P$.

Now for some $r\in M$ consider a haystack $H$ with respect to $r$. 
We want to add a ``needle'' to our haystack so that the resulting
subcomplex will again be a haystack. So we consider the set $A_r$ of
possible endpoints $v$ of edges through $r$, with  $v\prec r$:
\[
A_r \ \ :=\ \ \big\{v \in \Z^2\;\big|\; v \prec r 
    \mbox{ and  } H\cup\{[v,v + 2\vec{vr}]\}\mbox{ is a haystack }\big \}. 
\]
We want to prove  that $|A_r| \le 2$ for all $r\in M$.

\begin{figure}[ht]
\begin{center}
\input{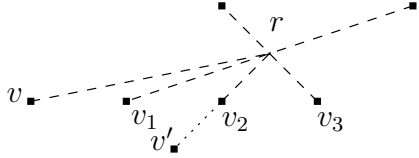}
\caption{\label{alpha.fig}Here $A_r =\{v_1,v_3\}$ and
  $A=\{v_1,v_2,v_3\}$, while by definition
 $v,v'\not\in A$ and $\alpha_1<\alpha_2<\alpha_3$.}
\end{center}
\end{figure}

We say that $v$ is \textbf{visible} from~$r$ if  the edge $[v,r]$
crosses no other edge or integral point.
Consider 
\[
A\ \ :=\ \ \big\{v\in \conv(\{r\}\cup A_r)\cap\Z^2\;\big|\;v 
\mbox{ is visible from }r\big\}.
\] 
As $v\in A_r$ is visible from $r$ we have $A \supseteq A_r$. 
Furthermore $ v\prec r$ holds for all   $ v\in A$. (See 
Figure~\ref{alpha.fig}).

We  now order $A$ by  the angles $\alpha(v)$ of $\vec{rv}$  with the
$x$-axis turning counter-clockwise and starting by $\pi$,  so that  we
have $\alpha_i=\alpha(v_i),\,\alpha_1 < \alpha_2 < \dots  <\alpha_k$. 
Indeed, we never have  $\alpha_i=\alpha_j$, otherwise
$r,v_i,v_j $ would lie on a line, but then one of the two points $v_i,v_j $
would not be visible from $r$, because both are $\prec r$.

Observe  that  $v_1\in A_r$: $ v\prec r$  for all   $ v\in A$, so a
point $ v$ with a smaller angle to the $x$-axis than the first one in
$A_r$  can't be in  $\conv(A_r\cup\{r\})\supset A$.

We say that a triangle  $[v_i,v_{i+1},r]$ is \textbf{empty} if there
is no edge through it and no $\frac{1}{2}\Z^2$-points in its interior.
The triangle $[v_i,v_{i+1},r]$ is empty 
 as $A$ contains all  points in
 $\conv(A_r\cup\{r\})\supset[v_i,v_{i+1},r] $ visible from $r$.
The midpoint  $s_i:=\frac{1}{2}(v_i+v_{i+1})$ is half-integer, 
$s_i\in M$, as $[v_i,v_{i+1},r]$ is empty. We also have $s_i\prec r$,
and so  $e_H(s_i) =[v_i,v_{i+1}]$, since the triangle $[v_i,v_{i+1},r]$ is
empty. Additionally it  has area 
$\frac{1}{4}$, otherwise the triangle wouldn't be empty.

Define $w_i:= r + \vec{v_ir}$ and 
$r':=\frac{1}{2}(v_1+w_2)$, $r'':=\frac{1}{2}(v_2+w_1) $. 
Then $v_1,w_2,v_2,w_1$ form a parallelogram with center $r$, and 
$r,r',r''$ are on a line (parallel to $(v_1v_2)$). 
So either $r'\prec r$ or $r'' \prec r$. 
\smallskip

\textbf{Case 1:} Suppose first that $r'\prec r$.\\
 The triangle $ \Delta=[v_1,v_2,w_2]$ is unimodular  as 
$\area(\Delta)= 2 \area[v_1,v_2,r]= \frac{1}{2} $; so there are no integer points between 
the  line  $(w_1w_2)$ and the line  $(v_1v_2)$. 
The edge  $e_H(r')$ has nonempty intersection with  these two lines 
(but doesn't cross $[v_1,w_1]$, since $v_1\in A_r$).
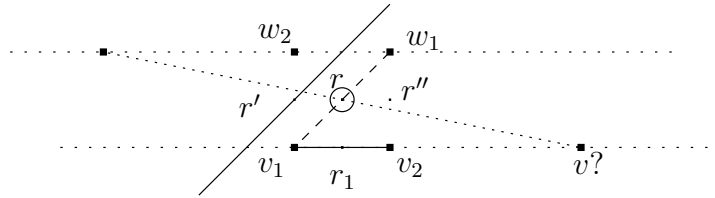
\begin{center}
\begin{picture}(0,0)%
\includegraphics{hay.unim-tri.pstex}%
\end{picture}%
\setlength{\unitlength}{3947sp}%
\begingroup\makeatletter\ifx\SetFigFont\undefined%
\gdef\SetFigFont#1#2#3#4#5{%
  \reset@font\fontsize{#1}{#2pt}%
  \fontfamily{#3}\fontseries{#4}\fontshape{#5}%
  \selectfont}%
\fi\endgroup%
\begin{picture}(4524,1224)(589,-1873)
\put(3046,-1711){\makebox(0,0)[lb]{\smash{\SetFigFont{12}{14.4}{\rmdefault}{\mddefault}{\updefault}{\color[rgb]{0,0,0}$v_2$}%
}}}
\put(2176,-1711){\makebox(0,0)[lb]{\smash{\SetFigFont{12}{14.4}{\rmdefault}{\mddefault}{\updefault}{\color[rgb]{0,0,0}$v_1$}%
}}}
\put(3106,-916){\makebox(0,0)[lb]{\smash{\SetFigFont{12}{14.4}{\rmdefault}{\mddefault}{\updefault}{\color[rgb]{0,0,0}$w_1$}%
}}}
\put(2176,-871){\makebox(0,0)[lb]{\smash{\SetFigFont{12}{14.4}{\rmdefault}{\mddefault}{\updefault}{\color[rgb]{0,0,0}$w_2$}%
}}}
\put(2626,-1786){\makebox(0,0)[lb]{\smash{\SetFigFont{12}{14.4}{\rmdefault}{\mddefault}{\updefault}{\color[rgb]{0,0,0}$r_1$}%
}}}
\put(4156,-1726){\makebox(0,0)[lb]{\smash{\SetFigFont{12}{14.4}{\rmdefault}{\mddefault}{\updefault}{\color[rgb]{0,0,0}$v ?$}%
}}}
\put(2626,-1186){\makebox(0,0)[lb]{\smash{\SetFigFont{12}{14.4}{\rmdefault}{\mddefault}{\updefault}{\color[rgb]{0,0,0}$r$}%
}}}
\put(3076,-1261){\makebox(0,0)[lb]{\smash{\SetFigFont{12}{14.4}{\rmdefault}{\mddefault}{\updefault}{\color[rgb]{0,0,0}$r''$}%
}}}
\put(2056,-1366){\makebox(0,0)[lb]{\smash{\SetFigFont{12}{14.4}{\rmdefault}{\mddefault}{\updefault}{\color[rgb]{0,0,0}$r'$}%
}}}
\end{picture}

\end{center}
But where could a third point $v \in  A_r $ (other than $v_1,v_3$) be? 
The line  $(r'r)$ is parallel to $(v_1v_2)$, we have 
$\alpha(r')<\alpha_1\leq\alpha_i$; and $r' \prec r, \; v\prec r$ for
all $v\in A_r$.
So all points of $A$ are on the same  side of $(r'r)$ as $v_1$ and $v_2$.   
So $v$  is on or beyond the  line $(v_1v_2)$ 
and  hence the edge through $r$ starting at $v$  would necessarily
cross the edge $e_H(r')$. So there can be no other point $v$  
in~$A_r$, that is, $|A_r|\leq 2$.
\smallskip

\textbf{Case 2:} The situation for $r''\prec r$ is similar:\\
The edge through $r''$ must be $e(r'')=[v_2,w_1]$, otherwise it would cut
$[v_1,w_1]$ or $[v_1,v_2]$; in the first case we would have $v_1\notin
A_r$ and in the second case  $v_2$ wouldn't be visible from~$r$.
And $[v_1,v_2,w_1]$ is again unimodular, 
so there is no possibility for a third $v\in A_r$.
\endofpr

Our Theorem 1 and its proof clearly extend to a more general
situation, namely the case of
a not necessarily simply connected lattice polygon (which may have
holes), possibly with additional, fixed inner edges. 

We can define the capacities $c_{m,n} := \frac{\log_2f(m,n)}{mn}$;
see~\cite{VoZi}.
From sublinearity of $f(m,n)$ it follows by Fekete's lemma
\cite[p.~85]{LW} that the limit capacities 
\[c_m\ \ :=\ \ \lim_{n\raw\infty}\frac{\log_2f(m,n)}{mn},\quad
  c_\Delta\ \ :=\ \ \lim_{n\raw\infty}\frac{\log_2f(n,n)}{n^2}
\]
exist.
 Theorem~\ref{thm:hay} yields the upper bounds
\[
c_m\ \ \le\ \ 3- \frac{1}{m},
\]
which includes the best known upper bounds
for all $c_m$ (compare \cite{VoZi}).

In generating triangulations with the ``haystack approach''
as in the proof of Theorem 1, one will in many situations 
have $|A_r|=1$. So probably our upper bound $c_\Delta\leq 3$ for the
limit capacity  $ c_\Delta$  is not sharp. 

As for lower bounds, the recursion formulas for
($n\times2$)- and ($n\times3$)-strips, as given in \cite{VoZi},
together with submultiplicativity, show  that 
$c_\Delta> 2.055$.


\thebibliography{9}\itemsep=0pt
\bibitem{Ore} {\sc S. Yu. Orevkov:} {Asymptotic number of triangulations with vertices in~${\mathbb Z}^2$,}
{{\sl J. Combinatorial Theory, Ser.~A} {\bf 86} (1999), 200-203.} 
\bibitem{VoZi} {\sc V. Kaibel \& G. M. Ziegler:} {Counting unimodular lattice triangulations,} Preprint, TU Berlin, November 2002, 31 pages;
  \url{arXiv:math.CO/0211268}; to appear in {\sl British Combinatorial Surveys}
  (C. D. Wensley, ed.), Cambridge University Press 2003.
\bibitem{LW}{\sc J. H. van Lint \& R. M. Wilson:}
{A Course in Combinatorics,}
{\sl Cambridge University Press}, Cambridge 1992.
\end{document}